
\documentstyle{amsppt}

\magnification=\magstep1

\hsize = 6.3 truein 
\vsize = 8.5 truein
\topskip = .25in

\document

\topmatter 
\title Bihomogeneity and Menger manifolds 
\endtitle
\author Krystyna Kuperberg
\endauthor  
\address  Department of Mathematics, Auburn University,
Auburn, AL 36849, USA 
\endaddress 
\thanks This research was supported in part by NSF grant
\# DMS-9401408
\endthanks
\email  kuperkm\@mail.auburn.edu
\endemail 

\abstract  It is shown that for every triple of integers $(\alpha,
\beta,\gamma)$ such that $\alpha \geq 1$, $\beta \geq 1$, and
$\gamma \geq 2$, there is a hom\-o\-geneous,
non-bihom\-o\-geneous  continuum whose every point has a
neighborhood homeomorphic  the Cartesian product of Menger
compacta $\mu^\alpha\times\mu^\beta\times\mu^\gamma$. In
particular, there is a homogeneous, non-bihom\-o\-geneous, Peano
continuum  of covering dimension four.  \endabstract 

\endtopmatter

\heading Introduction 
\endheading

A space is {\it $n$-hom\-o\-geneous\/} if for any pair of
$n$-point sets the  space admits a homeo\-morphism sending one
of the sets onto the other. A {\it hom\-o\-geneous\/} space is a
1-hom\-o\-geneous space. A space $X$ is {\it bihom\-o\-geneous\/}
if for any pair of points $p$ and $q$ in $X$,  there is a
homeo\-morphism $h: X\to X$ such that $h(p)=q$ and $h(q)=p$. A
space $X$ is {\it strongly locally hom\-o\-geneous\/} if for
every $p\in X$ and every neighborhood $U_p$  there is a
neighborhood $V_p$ such that for every  $q\in V_p$ there is a 
homeo\-morphism $h: X\to X$ with $h(p)=q$ and $h(x)=x$ for
$x\notin U_p$.   A continuum is a  compact, 
connected, metric  space containing more than one point. 

Around 1921, B.~Knaster  asked whether hom\-o\-geneity implies
bihom\-o\-geneity, and  C.~Kuratowski (Kazimierz Kuratowski)
\cite{Kur} gave an example of a   1-dimensional, non-locally
compact,  hom\-o\-geneous, non-bihom\-o\-geneous subset of the
plane. An example   similar to that of Kuratowski  can be
easily described as follows: Let $p$ and $q$ be distinct points
in the same composant of a nontrivial solenoid $\Sigma$.  The
composant of $\Sigma -\{p\}$ containing $q$ is hom\-o\-geneous
but it is not bihom\-o\-geneous; one end of the composant is
dense in it whereas  the other end is not, making swapping points
impossible.  It is also easy to obtain nonmetric examples. In
1986,  H.~Cook \cite{Cook} described a locally compact,
2-dimensional, hom\-o\-geneous, non-bihom\-o\-geneous metric
space. It is still not known if there is a 1-dimensional
[locally] compact metric example.

In 1930, D.~van~Dantzig \cite{Dan} restated Knaster's question
for continua, which was answered in \cite{KuK} by an example of
a locally connected, hom\-o\-geneous, non-bihom\-o\-geneous
continuum. The construction of the example starts with a  space
that is bihomogeneous but has a property  in some sense
contrary to bihom\-o\-geneity. This space  consists of
compatibly oriented circular fibers  such that any
homeomorphism  maps a fiber onto a fiber, and to swap certain
fibers the homeomorphism must reverse the  orientation of the
fibers. The next step is to replace each  fiber with a
hom\-o\-geneous space containing it as a retract and not
admitting a homeomorphism reversing the orientation of the
original fiber. The fibers are held together by a rigid grid
that is locally homeomorphic to the Cartesian product of two
Menger universal curves. Each circle is then replaced by a
larger fiber, a manifold, which contains $S^1$ as its retract,
but admits no homeomorphism changing the sign of the generator
of the first homology group represented by this $S^1$. The
dimension of the example equals  the dimension of the manifold 
plus two, giving a
 7-dimensional continuum. The seemingly unrelated property of
local connectedness is important for the notion of
2-hom\-o\-geneity: G.~S.~Ungar~\cite{Un} proved that
2-hom\-o\-geneous continua are locally connected. However,
hom\-o\-geneity does not imply 2-hom\-o\-geneity for Peano
continua \cite{KKT1, Ke1, Ke2, KKT2, Gar} and, as the above
example shows,  it does not imply bihom\-o\-geneity.

A substantially simpler, although not locally connected,
example of a hom\-o\-geneous, non-bihom\-o\-geneous continuum 
was given by P.~Minc  \cite {Minc}.  The ``model'' space of Minc's
example is a solenoid, whose each arc component is replaced by 
a sequence of ``glued together'' mapping cylinders of a degree
$m\geq 2$ map of $S^1$ onto $S^1$. For most pairs of composants
of a solenoid, a homeomorphism swapping the composants must be
orientation reversing (see \cite {Minc}). Hence in the above
continuum, not every two arc components can be swapped. To
achieve homogeneity, Minc takes the Cartesian product of this
continuum and the Hilbert  cube. To get a finite dimensional
example, a  manifold of the same homotopy type as the above
mapping cylinder can be used to replace the solenoid  composants.
Recently, K.~Kawamura \cite{Kaw} noticed  that using  Menger
manifolds and an ``$n$-ho\-mo\-topy  mapping cylinders'' (see
\cite{CKT and CKS}), the dimension of Minc's example can  be 
lowered to 2.

This paper shows that by applying Kawamura's idea  to  the 
construction of \cite{KuK},  a locally connected,
hom\-o\-geneous, non-bihom\-o\-geneous continuum  of dimension 4
can be obtained. The factorwise rigidity of the Cartesian
products of Menger compacta, immediately gives such examples in
all dimensions greater or equal to four.

\heading 1. Factorwise rigidity\endheading

K. Menger \cite{Men} defined $n$-dimensional universal compacta
in terms of the intersection of a sequence of polyhedra in
${\Bbb R}^{2n+1}$.  R.~D.~Anderson \cite{An1, An2} proved that
the 1-dimensional universal compactum, the Menger universal
curve $\mu ^1$, is homogeneous, and strongly locally
homogeneous. Not much was known about the higher dimensional
Menger universal   compacta until  M.~Bestvina \cite{Bes} 
characterized  the Menger universal compactum $\mu ^n$ as a space
that is topologically defined as follows:

\roster 

\item $\mu ^n$ is a compact n-dimensional metric space, 

\item $\mu ^n$ is  $LC^{n-1}$, 

\item $\mu ^n$ is  $C^{n-1}$,

\item $\mu ^n$ satisfies the Disjoint n-Disk Property, $DD^nP$.

\endroster

By \cite{Bes}, the compacta $\mu ^n$ are homogeneous  and
strongly locally homogeneous. An $n$-{\it dimensional Menger
manifold\/}, i.e., $\mu^n$-{\it manifold\/},  is a metric space
whose every point has a neighborhood homeomorphic to  $\mu^n$.

\definition{Definition} The Cartesian product
$X=\Pi_{\lambda\in\Lambda}X_\lambda$ is  {\it factorwise rigid\/}
if  every homeomorphism $h:X\to X$ preserves the Cartesian
factors: specifically, there is a permutation $\tau: \Lambda
\to \Lambda $ and  homeomorphisms $h_\lambda:X_{\tau(\lambda
)}\to X_\lambda$ such that if $h(\langle x_\lambda \rangle
)=\langle y_\lambda \rangle$, then $y_\lambda = h_\lambda
(x_{\tau (\lambda )})$. \enddefinition

\definition{Definition} An $i$-{\it fiber\/} of
$X=\Pi_{\lambda\in\Lambda}X_\lambda$ is a subset of $X$ of the form
$\Pi_{\lambda\in\Lambda}A_\lambda$, where $A_i=X_i$ and the
remaining factors are single points. An $i$-{\it cofiber\/} of
$X=\Pi_{\lambda\in\Lambda}X_\lambda$ is a subset of $X$ of the form
$\Pi_{\lambda\in\Lambda}A_\lambda$, where $A_i$ is a single
point and $A_j=X_j$ for $j\not= i$. \enddefinition

\definition{Definition} Points $x,y\in X$ are {\it
homologically separated in dimension $n$\/}, if they have
respective neighborhoods $U_x$ and $U_y$ such that   
$$i_*(\check H_n(U_x))\cap j_*(\check H_n(U_y))=0,$$ where $i:
U_x\hookrightarrow X$, $j: U_y\hookrightarrow X$ are the
inclusions, and $\check H$ is the $n$-th \v Cech homology group.
\enddefinition

The factorwise rigidity of the Cartesian product of two Menger
universal curves was first determined in \cite{KKT1} to show
that the product is not 2-homogeneous. The factorwise rigidity
of the Cartesian products of pseudo-arcs was showed in
\cite{BeKe} and \cite{BeLy}. J.~Kennedy Phelps
\cite{Ke1} proved that the Cartesian product of arbitrarily
many copies of $\mu^1$ is factorwise rigid, and by an
unpublished result of T.~Yagasaki, Kennedy's theorem extends to
the Cartesian products of copies of $\mu^n$ (see
\cite{CKT}, Section 3). D.~J.~Garity \cite{Gar} used the
K\"unneth and Eilenberg-Zilber formulas to show that finite
products of at least two  Menger universal compacta (of equal
or different dimensions, but excluding the product with all
factors $\mu^0$) are not 2-homogeneous. His proof is very close
to imply  factorwise rigidity. The notion of homology
separation was introduced in \cite{KKT2}. 

For dimensional reasons, any two points in $\mu^n$ are 
homologically separated in dimension $n$. At every $x\in
\mu^n$, there are arbitrarily small spheres $S^n$ embedded in
$\mu^n$ as retracts. Let $X=\mu^{n_1}\times \cdots \times
\mu^{n_k}$. The $m$-cycles, where $m=n_1 + \cdots +n_k=\dim
X$,  carried by two disjoint tori of form $S^{n_1}\times
\cdots\times S^{n_k}$ are not  homologous. Therefore, if
$h_1,h_2:X\to X$ are isotopic homeomorphisms, then $h_1=h_2$. 
The lemma below has analogs in the above mentioned
papers, but treats factorwise rigidity as a local property.

\proclaim{Lemma 1} Let  $X=X_1\times\cdots\times X_k$, where 
$X_i$ is homeomorphic to $\mu^{\alpha _i}$, $1\leq  \alpha
_1\leq  \cdots \leq   \alpha _k$. Let
$U=U_1\times\cdots\times U_k $,  $U_i\subset X_i$,   be an open
connected subset of $X$, and let  $h:U\to X$ be an  open
embedding.  Then  $h(x_1, \ldots, x_k)= (h_1(x_{\tau
(1)}),\ldots,h_k(x_{\tau (k)})) $, where $\tau$ is a
permutation and $h_i:U_{\tau (i)}\to X_i$ is an embedding.
\endproclaim

\demo{Proof}  Let $x=(x_1,\ldots , x_k)$ and $\bar x=(x_1,\bar
x_2 \ldots , \bar x_k)$ be two points in the same $1$-cofiber
of $U$. There is a sequence of $\alpha _1$-dimensional spheres
$C_n$ in $U_1$ containing $x_1$ with diam\/($C_n)\to 0$. Let
$K_n$ $=$ $\{ (y^n,x_2,\ldots ,x_k)\/ |\/ y^n\in C_n\}$ and
$\bar K_n$ $=$ $\{ (y^n,\bar x_2,\ldots ,\bar x_k)\/ |\/ y^n\in
C_n\}$.  The spheres $K_n$ are retracts of $X$, and the spheres 
$h(K_n)$ are retracts of $h(U)$. Every point of  $h(U)$ has small
neighborhoods $V_1\times \cdots \times V_k$  in $h(U)$,  which are
retracts of $X$.  Hence sufficiently small  spheres  $h(K_n)$ are
retracts of $X$. For some $i$ and infinitely many $n$'s,  $\pi
_i\circ h:K_n\to  X$ is essential.  Note that $\alpha _1=\alpha _i$,
so by the classical Hurewicz theorem, an essential map $S^{\alpha
_1}\to \mu^{\alpha _i}$ is homologically nontrivial. If the
$i$-coordinates of $h(x)$ and $h(\bar x)$ are different, then for
sufficiently large $n$, the images  $\pi _i(h(K_n))$ and  $\pi
_i(h(\bar K_n))$, are disjoint. Since the 1-cofibers of $U$ are
arcwise connected, there is a tube $S^{\alpha_1}\times [0,1]$  in $U$
joining $K_n$ and $\bar K_n$, which implies that the nontrivial
cycles represented by $\pi _i(h(K_n))$ and  $\pi _i(h(\bar K_n))$ are
homologous, contradicting the fact that distinct points in 
$\mu^{\alpha _i}$ are homologically separated in dimension ${\alpha
_i}$. Therefore the points $h(x)$ and $h(\bar x)$ are  in the same
$i$-cofiber of $X$.  Hence every 1-cofiber in $U$ is mapped into a
cofiber in $X$. By continuity, 1-cofibers close to a 1-cofiber  mapped
into an $i$-cofiber are mapped into  $i$-cofibers for the same $i$.
The compositions of the natural inclusions of $U_2\times \cdots \times
U_k$ into 1-cofibers, $h$ restricted to these 1-cofibers, and the
projection $\pi _i$ are isotopic, so they are identical. The proof
is completed by induction. \qed\enddemo

\proclaim{Corollary} \cite{See Problem 2, KKT2.}  Finite
Cartesian products of Menger universal compacta (also Menger
manifolds) are factorwise rigid.  \endproclaim

\heading 2. Grids\endheading

A homeomorphism $h:X\to X$ is {\it periodic\/} with {\it period\/}
$k\geq 1$ if $h^k(p)=p$ for every $p\in X$, but for every
$1\leq i<k$ and $p\in X$, $h^i(p)\not=p$. A closed subset $A$
of a compact metric space $X$ is a  $Z$-{\it set\/} if for every 
$\epsilon >0$,  there is a map $f:X\to X$  $\epsilon$-close to
the identity  with $f(X)\cap A = \emptyset$.

For positive integers $\alpha $, $\beta$, and $k\geq 2$, choose
$M$, $f_M$, $N$, $A$, $B$, $f_N$, $F$, and $Q$ as follows. 

\roster

\item $f_M:M\to M$, where $M$ is a $\mu^\alpha$-manifold, is
a periodic homeomorphism  of period $k$,

\item  $N=\mu^\beta$, 

\item $A$ and $B$ are disjoint, nonempty, homeomorphic $Z$-sets
in $N$, and $f_N:A\to B$ is a homeomorphism, such that the
quotient space of $N$ obtained by identifying each point $n\in A$
with its image $f_N(n)\in B$, is a $\mu^\beta$-manifold,  

\item  $F:M\times A\to M\times B$ is given by $F(m,n)=(f_M(m),
f_N(n))$,

\item $Q=(M\times N)/F$ is the quotient space obtained by  
identifying each point $(m,n)$ with the point $F(m,n)$.

\endroster

We refer to the continuum $Q$ as the $k$-{\it grid\/}. By a
slight abuse of notation, points in the quotient space are
denoted in the same way as the corresponding points in $M\times
(N-A)$. 

For $p=(m_p,n_p)\in Q$, let $$M_p=\{(m,n)\in Q \, |\,
n=n_p\},$$  $$N_p=\{(m,n)\in Q \,|\, m=f_M^i(m_p), i=0,\dots
,k-1\},$$  $$O_p=M_p\cap N_p.$$

\medskip

Call the sets $M_p$ and  $N_p$ {\it horizontal\/} and {\it
vertical fibers\/} respectively. The intersection of a horizontal
fiber and a vertical fiber is  a {\it necklace\/} and its
elements are {\it beads\/}. Note that the number of beads on each
necklace $O_p$ is $k$. 

\proclaim{Lemma 2} A homeomorphism $h:Q\to Q$ takes each
horizontal  and  vertical fiber  onto a horizontal or vertical
fiber, and a necklace  onto a necklace. \endproclaim

\demo{Proof} Using Lemma 1, the proof is identical to the proof
of Lemmas 5 and 6 in \cite{KuK}, where this is shown for the
case $\alpha =\beta =1$ and a specific $k$. \qed\enddemo

There is a cyclic order of a given set of beads on a necklace,
which cannot be arbitrarily disturbed by a homeomorphism of
$Q$. Let  $\phi :Q\to Q$ be the homeomorphism given by  $\phi
(m,n)=(f_M(m),n)$. For a point $p_0$ in $Q$,  denote by $p_i$
the point $\phi^i(p_0)$. Thus the necklace $O_{p_0}$ is the set
$\{p_0,\ldots ,p_{k-1}\}$.

\proclaim{Lemma 3} \cite{Compare with  Lemma 7 in KuK} Suppose
that  $h:Q\to Q$ is a homeomorphism and  $h(p_0)=p_0$. Then
there is an  $s$ such that $h(p_i)=p_{(si)\mod  k}$   for
$i=0,\ldots , k-1$. \endproclaim

\demo{Proof}  There is an arc $L_0$ joining $p_0$ and $p_1$
such that distinct $\phi^i(L_0)$ and $\phi^j(L_0)$ do not
intersect except for a possible common end point. Denote
$\phi^i(L_0)$ by $L_i$, and by $K$ the simple closed curve 
$\bigcup _{i=0}^{k-1} L_i$. Note that  $K$  is the union of
necklaces. Hence $h(K)$ is the union of necklaces; if $p\in
h(L_0)$ then $O_p\subset h(K)$. We have 
$$h(K)=\bigcup _{i=0}^{k-1}h( L_i)=\bigcup _{i=0}^{k-1}
\phi^ih(L_0).$$  The points $h(p_i)$ are ordered on the simple
closed curve $h(K)$ in such a way that the difference modulo $k$
in the indices between $h(p_i)$ and $h(p_{i+1})$) is a
constant. Therefore if $\phi(p_1)=p_s$, then 
$h(p_i)=p_{(si)\mod  k}$. \qed\enddemo

\medskip

\proclaim{Lemma 4} If  $h:Q\to Q$ is a homeomorphism and 
$h(p_0)\in O_{p_0}$, then there are integers $r$ and $s$  such that
$h(p_i)=p_{(r+si)\mod k}$   for $i=0,\ldots , k-1$. 
\endproclaim

\demo{Proof} Take $r$ and $s$ such that $h(p_0)=p_r$ and $\phi
^{-r}\circ h(p_1)=p_s$. \qed\enddemo

\heading 3. Circular fibers and  fiber replacing  \endheading

The next step is to construct a continuum built on the $k$-grid
$Q$  obtained from $Q\times I$ by the
identification $(p,1)=(\phi ^b(p),0)\/$, where $k$ is the
product of two positive integers $a$ and $b$.  Each of the sets
$O_p\times I$ transforms into $b$ circles called {\it circular
fibers\/}; $X$ decomposes into pairwise disjoint copies of $S^1$.
The number of beads of $O_p\times \{0\}$ on each circle $C$ is
$a$. The order of the beads determines the orientation of $C$. By
Lemma 3.1 of \cite{KKT2}, we have:

\proclaim{Lemma 5} A homeomorphism $h:X\to X$ takes each
circular fiber onto  a circular fiber. \endproclaim

Every point of $X$ has a neighborhood homeomorphic to the 
Cartesian product $M\times N\times I$. Orient the  $I$-fiber of
$Q\times I$ and transfer the orientation to the circulars
fibers of $X$. Similarly as in \cite{KuK} (see Lemma 11), a
homeomorphism of $X$ onto itself either preserves  orientation on
all circular fibers (it is then {\it orientation preserving\/}),
or reverses  orientation on all circular fibers (it is then {\it
orientation reversing\/}). 

\medskip

 As in the previous section let $O_{p_0}=\{p_0,\ldots
,p_{k-1}\}$ be a necklace in $Q$. Denote the point $(p_i,0)\in
X$ by $q_i$.

\proclaim{Lemma 6} If $a\geq 3$, then  every homeomorphism
$h:X\to X$ such that $h(q_0)=q_1$, $h(q_1)=q_0$,  and
$h(Q\times \{0\})=Q\times \{0\}$ is orientation reversing.
\endproclaim

\demo{Proof} By Lemma 4, $h(q_i)=q_{(1-i)\mod ab}$. Suppose
that $h$ is orientation preserving. Then
$h(\{q_0,q_b,q_{2b},\ldots \})$ $=$
$\{q_1,q_{1+b},q_{1+2b},\ldots \}$ preserving order. Hence
$h(q_b)=q_{1+b}$ and $1+b=(1-b)\mod ab$. So $2 \mod a =0$, which
is a contradiction. \qed\enddemo

\medskip

If one were to follow the procedure described in \cite{KuK}, each
circular fiber $C$ of $X$ would be replaced by a manifold
$E$ which contains $C$ as its retract, and such that every
autohomeomorphism of $E$  takes the element  of the first homology
group  represented by $C$ onto itself; in particular it does not
change the sign of this element. (Note that in \cite{KuK}, $C$, $E$,
and $X$  are denoted by different symbols.)  The resulting continuum
$D$ is the union of pairwise disjoint copies  of the same manifold
$E$, called manifold fibers. Since in in \cite{KuK} $a=b=3$; each
circular fiber in $X$ consists of three segments; each manifold fiber
of $D$ consists of three identical pieces.  $D$
contains a copy of $X$ as its retract. Although $X$ needs not be
invariant under a homeomorphism $h:D\to D$, $h$  induces a
homeomorphism of $X$ preserving the correspondence of the circular
fibers to the manifold fibers given by the inclusion
$X\hookrightarrow D$. It is shown that:

\roster  

\item every homeomorphism $h:D\to D$ maps a  manifold  fiber
onto a  manifold fiber, 

\item  some  manifold fibers  cannot be swapped by a
homeomorphism of $D$. 

\endroster

Kawamura's idea \cite{Kaw} to modify Minc's example \cite{Minc}
can be also applied to modify $D$. Instead of the manifold $E$
take a Menger manifold $\Omega $ consisting of $a$ identical pieces
homeomorphic to a $\mu ^\gamma$-manifold $P$, where $\gamma
\geq 2$. $P$ corresponds to the  mapping cylinder of a degree
two map of $S^1$ onto $S^1$.  The Menger manifold $\Omega $
is similar to the Menger manifold $L_n$ in \cite{Kaw}, Section
3,  and has the following properties:

\roster

\item  $\Omega =\bigcup _{i=0}^{a-1} \Omega _i$, and there are 
homeomorphisms $\tau_i :P\to\Omega _i$. 

\item  $\Omega _i\cap \Omega _j= \emptyset$ if $|(i-j) \mod a|\not= 1$.

\item  There are two disjoint homeomorphic $Z$-sets in $P$,
$P_0$ and $P_1$, such that $\Omega _i\cap \Omega _j= \tau_i (P_0)=\tau_j
(P_1) $ if $(i-j) \mod a= 1$, and  $\Omega _i\cap \Omega _j= \tau_i
(P_1)=\tau_j (P_0) $ if $(j-i) \mod a= 1$.

\item There is a simple closed curve $K\subset\Omega $ intersecting
each $\Omega _i$ in an arc such that $K$ is a retract of $\Omega $, every
homeomorphism  $\Omega \to\Omega $  takes the element  of the first
homology group  represented by (oriented) $K$ onto itself without
changing its sign.

\endroster

The Menger manifold $P$ can be used to replace the circular
fiber of the continuum $X$ to obtain a non-bihomogeneous 
continuum $Y$. Namely, $Y$ is the quotient space obtained from
$Q\times P$ by identifying each point $(p,x)$ with $(\phi ^b
(p), \tau^{-1}_{(i+1)\mod a}\circ\tau_i(x))$, where $Q$ is the
$k$-grid considered in the beginning of this section, $p\in Q$,
and $x\in P_1$. Note that $\tau^{-1}_{(i+1)\mod a}\circ\tau_i
:P_1\to P_0$.

\proclaim{Lemma 7} The continuum $Y$ is homogeneous.
\endproclaim

\demo{Proof} The proof uses strong local homogeneity of the Menger
compacta and is similar the proofs of Lemmas 3 and 4 in \cite{KuK}.
\qed\enddemo

\proclaim{Lemma 8} If $a\geq 3$ and $b\geq 2$, then $Y$ is not
bihomogeneous. \endproclaim

\demo{Proof} Since $Y$ is locally homeomorphic to $\mu^\alpha
\times\mu^\beta \times \mu^\gamma$, the local factorwise
rigidity holds. Using the local Cartesian product structure, we
may define the $\mu^\alpha$-, $\mu^\beta$-, and
$\mu^\gamma$-fibers, as well as the  $\mu^\alpha$-,
$\mu^\beta$-, and $\mu^\gamma$-cofibers. The 
$\mu^\gamma$-fibers are homeomorphic to $P$, and the
$\mu^\gamma$-cofibers are homeomorphic to $Q$.  The necklaces
in the $\mu^\gamma$-cofibers have $ab$ beads, whereas the 
necklaces in the $\mu^\alpha$- and $\mu^\beta$-cofibers contain
only $a$ beads. Therefore every homeomorphism $Y\to Y$ maps the
$\mu ^\gamma$-fibers onto $\mu ^\gamma$-fibers and $\mu
^\gamma$-cofibers onto $\mu ^\gamma$-cofibers even if $\alpha
=\beta =\gamma$. Since $b \geq 2$, there are at least two
distinct $\mu^\gamma$-fibers passing through the same necklace
of $Q\times \{x\}$, where $x$ is a point in $P-P_1$. By Lemma 6 and
by the above consideration, there are two such fibers that cannot be
swapped.
\qed\enddemo

The above lemma could be compared to Lemmas 13,
14, 15, and 16 in \cite{KuK}.  However, factorwise rigidity
involving all three factors,
$\mu^\alpha$, $\mu^\beta$, and $\mu^\gamma$, makes the proof much
simpler. 

\proclaim{Theorem}  For every triple of integers $(\alpha,
\beta,\gamma)$ such that $\alpha \geq 1$, $\beta \geq 1$, and
$\gamma \geq 2$, there is a hom\-o\-geneous,
non-bihom\-o\-geneous  continuum whose every point has a
neighborhood homeomorphic to 
$\mu^\alpha\times\mu^\beta\times\mu^\gamma$.  \endproclaim

\remark{Remark} $Y$ is the quotient space of $\mu^\alpha
\times\mu^\beta \times \mu^\gamma$ with some identifications
made along $Z$-sets. It is not known whether the same effect
can be achieved on the product $\mu^1 \times\mu^1 \times
\mu^1$. \endremark

\definition{Question 1} Does there exist  a hom\-o\-geneous,
non-bihom\-o\-geneous, Peano continuum of dimension lower than
4?  \enddefinition

\definition{Question 2} Does there exist  a hom\-o\-geneous,
non-bihom\-o\-geneous, Peano continuum whose every point has a
neighborhood homeomorphic to  $\mu^1 \times\mu^1
\times \mu^1$ ? \enddefinition

\remark{Remark} While the second question remains open, since the
submission of this paper, the first question has been answered by G.
Kuperberg  \cite{KuG}.  For any pair of integers
$\alpha,\beta $ such that $\alpha \geq 1$ and $\beta \geq 2$, he
constructs a hom\-o\-geneous, non-bihom\-o\-geneous Peano continuum
with the local structure of  $\mu^\alpha\times\mu^\beta$.
\endremark

\enddocument